# Rényi entropy for multivariate controlled autoregressive moving average systems[1]

Salah H. Abid, Uday J. Quaez, Javier E. Contreras-Reyescor

## Abstract

Rényi entropy is an important measure in the context of information theory as a generalization of Shannon entropy. This information measure was often used for uncertainty quantification of dynamical behaviour of stochastic processes. In this paper, we study in detail this measure for multivariate controlled autoregressive moving average (MCARMA) systems. The characteristic function of output process is represented from the terms of its residual characteristic function. An explicit formula to compute the Rényi entropy for the output process of MCARMA system is derived. In addition, we investigate the covariance matrix to find the upper bound of Rényi entropy. We present three simulations that serve to illustrate the behavior of information in MCARMA system, where the control and noise follow the Gaussian, Cauchy and Laplace distributions. Finally, the behaviour of Rényi entropy is illustrated in two real-world applications: a paper-making process and an electric circuit system.

## 1 Introduction

Autoregressive moving average (ARMA) processes [Brockwell and Davis(1991)] are strong tools for understanding, analyzing, representing and predicting future values of a very wide range of phenomena. [Andel(1983)] was one of the first authors whose were interested by determining the marginal distribution of univariate $\mathrm{AR}(1)$ process $x_t = \phi x_{t-1} + w_t$. He derived the characteristic function of this process in the following form

$$\varphi_{x_t}(s) = \prod_{j=0}^{\infty} \varphi_{w_t}(\phi^j s),$$

where $\varphi_{x_t}(s)$ and $\varphi_{w_t}(s)$ are the characteristic functions of $x_t$ and $w_t$, respectively. [Priestley(1983)] represented the characteristic function of univariate $\mathrm{AR}(2)$ by its residual characteristic functions. The behavior of univariate $\mathrm{AR}(1)$ process is studied by [Simy(1986)], where the Gamma and Weibull densities were used as perturbations. [Abid(2006)] derived a general formula to represent the characteristic function of univariate $\mathrm{ARMA}(p, q)$ process.

Parameter estimation technology based on input and output data has various applications in circuit system modeling, model prediction, automatic control, and system identification

System identification is an important issue to model dynamical/complex systems based on input and output data, such as in circuit system modeling, model prediction, and automatic control [Bao et al.(2011), Li et al.(2023)]. The controlled ARMA system is one of the most considered approach for system identification, which played a major role in many industrial applications such as speech processing [Makhoul(1975), Rabiner and Schafer(2005), Markel and Gray(2013)], biological signal processing, electroencephalogram classification EEG [Gersch et al.(1980)], electrocardiogram classification ECG [Vuksanovic and Alhamdi(2013)], geophysical exploration [Robinson et al.(1986), Robinson and Treitel(2000)], communication systems [Ranade and Sahai(2015)], and in many other applications. For example, [Sadeghi et al.(2017)] considered an ARMAX dynamic model (an univariate version of multivariate CARMA one) for a solar power plant to describe its characteristics, [Liu et al.(2018)] considered a multivariate stochastic system with autoregressive white noise for paper-making process modelled by multivariate CARMA (MCARMA) model [Bao et al.(2011)], [Hag ElAmin(2020)] also consider an MCARMA model for an electric circuit system with two-input voltage and a single-output voltage measured from a capacitor for a posterior clustering analysis, and more recently [Li et al.(2023)] considered a servomechanism system based on a magnet synchronous motor modelled as in latter works.

---

[1]

On the other hand, Shannon entropy [Cover and Thomas(2006)] is an information measure to quantify the uncertainty of an event. [Ràyi(1961)] generalized this measure for probability distribution which means the sensitive to the fine details of a density function [Contreras-Reyes(2022b)]. Recent important work done by [Abid and Hameed(2014)] considered the problem of deriving the Shannon entropy for univariate ARMA process. The main idea consists in find the expression of the characteristic function of output process in terms of its residual characteristic function, and then solve the problem of writing causal function parameters in terms of corresponding univariate ARMA$(p,q)$ process. More recently, [Gibson(2018)] studied the Shannon entropy power and mutual information of univariate AR$(p)$ process. [Contreras-Reyes(2022a)] computed the Rényi and differential entropies and Rényi and Kullback–Leibler divergences for a vector ARFIMA process using the characteristic and impulse response function. Regarding nonlinear processes, [Contreras-Reyes(2023)], [Contreras-Reyes(2024a)] and [Contreras-Reyes(2024b)] computed the differential entropy, the Kullback–Leibler and Jeffrey's divergences, and Jensen-variance distance for univariate and multivariate self-exciting threshold autoregressive and ARFIMA processes.

Parameter estimation techniques of MCARMA systems has been the main focus by researchers, where least squares estimation has been the most used iterative method for two-input single-output models [see Li et al.(2023), and references therein]. As was mentioned latter, extensive literature have been developed to quantify uncertainty in signal processing and time series modelling. Beyond measuring uncertainty, the information measures played an important role in summarize the overall information provided by estimated parameter of an specific model. For system identification based on input and output data, there is a lack in uncertainty and information quantification issues, where MCARMA system considered a large number of parameters related to polynomials orders. In addition, information measures also could help for system identification in the evaluation analysis of model performance [Contreras-Reyes(2024a)]. For latter reasons, in this paper we discuss the information in a MCARMA process according to Rényi measure. Given that covariance matrix of these processes can be calculated, the explicit expression of Rényi entropy of this model is derived using some conditions such as stationary, stability and independence of process. In addition, the formula of the characteristic function of output process from its residual characteristic functions is found.

Then, the upper bound of Rényi entropy is derived.

This paper is organized as follows. In Section 2, we begin with a preliminary material of the Rényi measure of information. Section 3 provides a description of MCARMA systems. Asymptotic expression for the characteristic function of output process of these systems is derived in Section 4. Then, we give the expression and upper bound of Rényi entropy in Sections 5 and 6, respectively. Numerical examples to illustrate this study are given in Section 7. Finally, the conclusion of this paper is shown in Section 8.

# 2 Preliminary

In this section, we introduce some basic definitions, notations and lemmas related to entropy theory. Some notations are introduced for convenience: $\mathbf{I}_d$ represents an $d \times d$ identity matrix, $N_d$, $C_d$ and $L_d$ represent multivariate Gaussian, Cauchy and Laplace distributions, respectively, in the space $\mathbb{R}^d$ and $\mathbf{A}^\top$ is the transpose of the matrix $\mathbf{A}$.

**Definition 1** *[Cover and Thomas(2006)] Let $\boldsymbol{X}$ be a continuous random vector in $\mathbb{R}^d$ with probability density function $f(x;\theta)$. Then the Shannon entropy is defined as*

$$H(\mathbf{X}) = -\mathbb{E}[\log\{f(\mathbf{x};\theta)\}], \tag{1}$$

where $\mathbb{E}$ and $\log$ denote expectation and natural logarithmic, respectively.

**Definition 2** *[Ràyi(1961)] An $\alpha$th-order Rényi entropy of a continuous random vector $\boldsymbol{X} \in \mathbb{R}^d$ with probability density function $f(\boldsymbol{x};\theta)$ is defined as*

$$R_\alpha(\mathbf{X}) = \begin{cases} \frac{1}{1-\alpha}\log\{\mathbb{E}[f(\mathbf{x};\theta)^{\alpha-1}]\}, & if \quad 0 < \alpha < \infty, \alpha \neq 1; \\ H(\mathbf{X}), & if \quad \alpha \to 1. \end{cases} \tag{2}$$

The following Lemma allows to obtain an upper bound for Rényi entropy for multivariate random variables, previously considered in [Contreras-Reyes and Cortà(2016)] and [Abid et al.(2021)].

**Lemma 1** *[S chez-Moreno et al.(2011)] If a random vector $\boldsymbol{X} \in \mathbb{R}^d$ has covariance matrix $\boldsymbol{S} = \mathbb{E}[\boldsymbol{X}\boldsymbol{X}^\top]$, then the upper bound of Rényi entropy of $\boldsymbol{X}$ in terms of $\boldsymbol{S}$ is*

$$R_\alpha(\mathbf{X}) \leq \mathcal{C}_d(\alpha) + \frac{1}{2}\log\{\mathrm{d}et(\mathbf{S})\}, \tag{3}$$

where

$$\mathcal{C}_d(\alpha) = \begin{cases} \frac{d}{2}\log\left(\frac{\pi b}{\alpha-1}\right) + \frac{1}{\alpha-1}\log\left(\frac{b}{2\alpha}\right) + \log\Gamma\left(\frac{\alpha}{\alpha-1}\right) - \log\Gamma\left(\frac{b}{2(\alpha-1)}\right), & if \ \alpha > 1\,; \\ \frac{d}{2}\log\left(\frac{\pi b}{1-\alpha}\right) + \frac{\alpha}{\alpha-1}\log\left(\frac{b}{2\alpha}\right) - \log\Gamma\left(\frac{\alpha}{1-\alpha}\right) + \log\Gamma\left(\frac{b}{2(1-\alpha)}\right), & if \ \alpha \in \left(\frac{d}{d+2}, 1\right)\,, \\ \frac{d}{2}\log(2\pi e), & if \ \alpha \to 1\,; \end{cases}$$

$b = (2+d)\alpha - d$ and $\Gamma(z) = \int_0^\infty t^{z-1}e^{-t}dt,\ z > 0.$

Depending of multivariate distribution of $\mathbf{X}$, Rényi entropy is computed with respect to order $\alpha$ as follows.

**Example 1** *[Song(2001), Contreras-Reyes(2015)] If $\boldsymbol{X} \sim N_d(\boldsymbol{\mu}, \boldsymbol{\Sigma})$ for some positive definite matrix $\boldsymbol{\Sigma}$, then the Rényi entropy of $\boldsymbol{X}$ can be written as*

$$R_\alpha(\mathbf{X}) = \begin{cases} \frac{1}{2}\log\{\det(2\pi\mathbf{S})\} - \frac{d}{2(1-\alpha)}\log(\alpha), & if \quad 0 < \alpha < \infty, \alpha \neq 1; \\ \frac{1}{2}\log\{\det(2\pi e\mathbf{S})\}, & if \quad \alpha = 1. \end{cases} \tag{4}$$

**Example 2** *[Abid and Quaez(2019)] If $\boldsymbol{X} \sim C_d(\boldsymbol{\mu}, \boldsymbol{\Sigma})$, then the Rényi entropy of $\boldsymbol{X}$ is*

$$R_\alpha(\mathbf{X}) = \begin{cases} \log\left(\frac{\alpha^{-\frac{d}{1-\alpha}}\sqrt{\det(4\pi\mathbf{S})}\Gamma(\frac{d+1}{2})}{\sqrt{\pi}}\right), & if \quad 0 < \alpha < \infty, \alpha \neq 1; \\ \log\left(\frac{\sqrt{\det(4e^2\pi\mathbf{S})}\Gamma(\frac{d+1}{2})}{\sqrt{\pi}}\right), & if \quad \alpha \to 1. \end{cases} \tag{5}$$

## 3 Multivariate controlled autoregressive and moving average systems

The MCARMA system has been introduced by [Bao et al.(2011)], which it can be mathematically described by a stochastic control system such as

$$\begin{cases} \mathbf{x}_t = \sum_{i=1}^{p} \mathbf{A}_i\mathbf{x}_{t-i} + \sum_{j=0}^{r} \mathbf{B}_j\mathbf{u}_{t-j} + \sum_{k=0}^{q} \mathbf{D}_k\mathbf{w}_{t-k}, \\ \mathbf{x}_0 = 0, \end{cases} \tag{6}$$

where $\mathbf{x}_t$ is the system output state, $\mathbf{u}_t$ is the system input vector, $\mathbf{w}_t$ is the white noise vector with mean zero. $\mathbf{A}_i \in \mathbb{R}^{d\times d}$, $\mathbf{B}_j \in \mathbb{R}^{d\times d}$, $\mathbf{D}_k \in \mathbb{R}^{d\times d}$, and $\mathbf{B}_0 = \mathbf{D}_0 = \mathbf{I}_d$ ($\mathbf{I}_d$ denotes identity matrix of dimension $d$). Assume that

$$\mathbb{E}[\mathbf{U}_t\mathbf{U}_s^\top] = \begin{cases} \mathbf{S}_\mathbf{u}, & if \quad s = t, \\ \mathbf{0}, & if \quad s \neq t; \end{cases} \tag{7}$$

$$\mathbb{E}[\mathbf{W}_t\mathbf{W}_s^\top] = \begin{cases} \mathbf{S}_\mathbf{w}, & if \quad s = t, \\ \mathbf{0}, & if \quad s \neq t; \end{cases} \tag{8}$$

and both processes $\mathbf{u}_t$ and $\mathbf{w}_t$ are independent for any $s$ and $t$. Rewriting the equation (6) by using the lag operator $L^j\mathbf{x}_t = \mathbf{x}_{t-j}$, we get

$$[\mathbf{I}_d - \mathbf{A}(L)]\mathbf{x}_t = \mathbf{B}(L)\mathbf{u}_t + \mathbf{D}(L)\mathbf{w}_t, \tag{9}$$

where $\mathbf{A}(L) = \sum_{i=1}^{p} \mathbf{A}_iL^i$, $\mathbf{B}(L) = \mathbf{I}_d + \sum_{j=1}^{r} \mathbf{B}_jL^j$, and $\mathbf{D}(L) = \mathbf{I}_d + \sum_{k=1}^{q} \mathbf{D}_jL^j$.

**Definition 3** *[Madsen(2007)] The output process $\boldsymbol{x}_t$ of the control system in (6) is said to be stable if and only if $det([\boldsymbol{I} - \boldsymbol{A}(z)]^{-1}) \neq 0$, for all complex numbers $z$, $|z| < 1$.*

Assume that the process $\mathbf{x}_t$ is stable, then the equation (9) can be written as

$$\mathbf{x}_t = \mathbf{M}(L)\mathbf{u}_t + \mathbf{M}^*(L)\mathbf{w}_t, \tag{10}$$

with notations $\mathbf{M}(L) = [\mathbf{I}_d - \mathbf{A}(L)]^{-1}\mathbf{B}(L) = \sum_{j=0}^{\infty} \mathbf{M}_j L^j$ and $\mathbf{M}^*(L) = [\mathbf{I}_d - \mathbf{A}(L)]^{-1}\mathbf{D}(L) = \sum_{k=0}^{\infty} \mathbf{M}_k^* L^k$. Then the equation (10) can be written as

$$\mathbf{x}_t = \sum_{j=0}^{\infty} \mathbf{M}_j L^j \mathbf{u}_t + \sum_{k=0}^{\infty} \mathbf{M}_k^* L^k \mathbf{w}_t. \tag{11}$$

Now, we compute the matrices $\mathbf{M}_j$, $j = 0,1, \ldots$, $\mathbf{M}_k^*$, $k = 0,1, \ldots$. By setting $\mathbf{M}_0 = \mathbf{M}_0^* = \mathbf{I}_d$ and multiplying both sides of the equation (11) by the operator $[\mathbf{I}_d - \mathbf{A}(L)]$, we get

$$\begin{aligned}[\mathbf{I}_d - \mathbf{A}(L)]\mathbf{x}_t &= \left[\mathbf{I}_d + \sum_{i=1}^{\infty} \left(\mathbf{M}_i - \sum_{j=1}^{i} \mathbf{A}_j \mathbf{M}_{i-j}\right) L^i\right]\mathbf{u}_t \\ &+ \left[\mathbf{I}_d + \sum_{i=1}^{\infty} \left(\mathbf{M}_i^* - \sum_{j=1}^{i} \mathbf{A}_j \mathbf{M}_{i-j}^*\right) L^i\right]\mathbf{w}_t. \end{aligned} \tag{12}$$

If $\mathbf{A}_j = \mathbf{0}$ for $j > p$, $\mathbf{B}_i = \mathbf{0}$ for $i > r$, $\mathbf{D}_i = \mathbf{0}$ for $i > q$, then the equation (6) can be written as

$$[\mathbf{I}_d - \mathbf{A}(L)]\mathbf{x}_t = \sum_{j=0}^{\infty} \mathbf{B}_j L^j \mathbf{u}_t + \sum_{k=0}^{\infty} \mathbf{D}_k L^k \mathbf{w}_t. \tag{13}$$

Finally, by comparing the coefficients in equations (12) and (13), we obtain

$$\mathbf{M}_i = \mathbf{B}_i + \sum_{j=1}^{i} \mathbf{A}_j \mathbf{M}_{i-j}, \quad i = 1,2, \ldots, r; \tag{14}$$

$$\mathbf{M}_k^* = \mathbf{D}_k + \sum_{j=1}^{k} \mathbf{A}_j \mathbf{M}_{k-j}^*, \quad k = 1,2, \ldots, q. \tag{15}$$

# 4 Characteristic function

Based on the equivalent equation (11) of MCARMA process $\{\mathbf{x}_t\}$, the representation of characteristic function from its residual is

$$\begin{aligned}\boldsymbol{\varphi}_{\mathbf{x}_t}(\mathbf{s}) &= \mathbb{E}\left[\exp\left(i\mathbf{s}^\top\left\{\sum_{j=0}^{\infty} \mathbf{M}_j \mathbf{u}_{t-j} + \sum_{k=0}^{\infty} \mathbf{M}_k^* \mathbf{w}_{t-k}\right\}\right)\right] \\ &= \prod_{j=0}^{\infty} \mathbb{E}\left[\exp(i\mathbf{s}^\top \mathbf{M}_j \mathbf{u}_{t-j})\right] \prod_{k=0}^{\infty} \mathbb{E}[\exp(i\mathbf{s}^\top \mathbf{M}_k^* \mathbf{w}_{t-k})], \end{aligned} \tag{16}$$

with $i = \sqrt{-1}$. Equation (16) is obtained using independence condition of $\{\mathbf{u}_t\}$ and $\{\mathbf{w}_t\}$ processes. Then,

$$\boldsymbol{\varphi}_{\mathbf{x}_t}(\mathbf{s}) = \prod_{j=0}^{\infty} \boldsymbol{\varphi}_{\mathbf{M}_j \mathbf{u}_{t-j}}(\mathbf{s}) \boldsymbol{\varphi}_{\mathbf{M}_j^* \mathbf{w}_{t-j}}(\mathbf{s}), \tag{17}$$

where $\boldsymbol{\varphi}_{\mathbf{M}_j \mathbf{u}_{t-j}}(\mathbf{s})$ and $\boldsymbol{\varphi}_{\mathbf{M}_j^* \mathbf{w}_{t-j}}(\mathbf{s})$ represent the characteristic functions of $\mathbf{M}_j \mathbf{u}_{t-j}$ and $\mathbf{M}_j^* \mathbf{w}_{t-j}$, respectively.

# 5 Rényi Entropy of MCARMA process

In this section, we drive the expression of Rényi entropy of MCARMA system based on multivariate Gaussian, multivariate Cauchy and multivariate Laplace distributions.

## 5.1 Multivariate Gaussian distribution

Assume that control and white noise processes are both multivariate Gaussian distributed with zero mean with covariance matrices $\mathbf{S}_\mathbf{u}$ and $\mathbf{S}_\mathbf{w}$, respectively; and they satisfy the independence conditions. Then from the equation (17) the characteristic function of process $\mathbf{x}_t$ can be written as

$$\boldsymbol{\varphi}_{\mathbf{x}_t}(\mathbf{s}) = \exp(-\tfrac{i}{2}\mathbf{s}^\top \sum_{j=0}^{\infty} (\mathbf{M}_j \mathbf{S}_\mathbf{u} \mathbf{M}_j^\top + \mathbf{M}_j^* \mathbf{S}_\mathbf{w} {\mathbf{M}^*}_j^\top)\mathbf{s}). \tag{18}$$

It can be seen from (18) that the distribution of process $\{\mathbf{x}_t\}$ is multivariate Gaussian with zero mean and covariance matrix

$$\boldsymbol{\Sigma} = \sum_{j=0}^{\infty} (\mathbf{M}_j \mathbf{S}_\mathbf{u} \mathbf{M}_j^\top + \mathbf{M}_j^* \mathbf{S}_\mathbf{w} {\mathbf{M}^*}_j^\top). \tag{19}$$

Considering Example 1, the Rényi entropy of process $\{\mathbf{x}_t\}$ based on system (6) is

$$R_\alpha(\mathbf{Y}) = \begin{cases} \frac{1}{2}\log[\det(2\pi\boldsymbol{\Sigma})] - \frac{d\log(\alpha)}{2(1-\alpha)}, & if \quad \alpha \neq 1; \\ \frac{1}{2}\log[\det(2\pi e\boldsymbol{\Sigma})], & if \quad \alpha \to 1. \end{cases} \tag{20}$$

## 5.2 Multivariate Cauchy distribution

Assume that control and white noise processes are both multivariate Cauchy distributed zero mean with scale matrices $\mathbf{S}_\mathbf{u}$ and $\mathbf{S}_\mathbf{w}$, respectively; and they satisfy the independence conditions. The characteristic function of any stationary process $\{\mathbf{x}_t\}$ with multivariate Cauchy distribution with location vector parameter $\boldsymbol{\mu}$ and scale matrix parameter $\boldsymbol{\Sigma}$ is [Bian(1989), Brockwell and Davis(1991)]

$$\boldsymbol{\varphi}_{\mathbf{x}_t}(\mathbf{s}) = \exp(i\mathbf{s}^\top\boldsymbol{\mu} - \sqrt{\mathbf{s}^\top\boldsymbol{\Sigma}\mathbf{s}}). \tag{21}$$

Therefore, from equation (17) the characteristic function of output process $\mathbf{x}_t$ based on control system (6) can be written as

$$\boldsymbol{\varphi}_{\mathbf{x}_t}(\mathbf{s}) = \exp\left(-\sum_{j=0}^{\infty}\left[\sqrt{\mathbf{s}^\top\mathbf{M}_j\mathbf{S}_\mathbf{u}\mathbf{M}_j^\top\mathbf{s}} + \sqrt{\mathbf{s}^\top\mathbf{M}_j^*\mathbf{S}_\mathbf{w}{\mathbf{M}_j^*}^\top\mathbf{s}}\right]\right). \tag{22}$$

Given that both matrices $\mathbf{S}_\mathbf{u}$ and $\mathbf{S}_\mathbf{w}$ are positive definite, then they can be written as $\mathbf{S}_\mathbf{u} = \mathbf{L}_\mathbf{u}\mathbf{L}_\mathbf{u}^\top$ and $\mathbf{S}_\mathbf{w} = \mathbf{L}_\mathbf{w}\mathbf{L}_\mathbf{w}^\top$, respectively; where $\mathbf{L}_\mathbf{u}, \mathbf{L}_\mathbf{w} \in \mathbb{R}^{d\times d}$. Therefore, the characteristic function of process $\mathbf{x}_t$ is

$$\boldsymbol{\varphi}_{\mathbf{x}_t}(\mathbf{s}) = \exp\left(-\sum_{j=0}^{\infty}\left[\sqrt{\mathbf{s}^\top\mathbf{K}_j\mathbf{K}_j^\top\mathbf{s}} + \sqrt{\mathbf{s}^\top\mathbf{K}_j^*{\mathbf{K}^*}_j^\top\mathbf{s}}\right]\right), \tag{23}$$

where $\mathbf{K}_j = \mathbf{M}_j\mathbf{L}_\mathbf{u}$ and $\mathbf{K}_j^* = \mathbf{M}_j^*\mathbf{L}_\mathbf{w}$, $j = 0,1,2,\ldots$.

Assume that matrices $\mathbf{K}_j\mathbf{K}_j^\top$ and $\mathbf{K}_j^*{\mathbf{K}^*}_j^\top$, $j = 0,1,2,\ldots$, are proportional one to other then there exists a matrix $\mathbf{C}$ such that [Bian(1989)]

$$\sqrt{\mathbf{s}^\top\mathbf{C}\mathbf{s}} = \sqrt{\mathbf{s}^\top\mathbf{K}_j\mathbf{K}_j^\top\mathbf{s}} + \sqrt{\mathbf{s}^\top\mathbf{K}_j^*{\mathbf{K}^*}_j^\top\mathbf{s}}. \tag{24}$$

Latter assumption gives the following result

$$\boldsymbol{\varphi}_{\mathbf{x}_t}(\mathbf{s}) = \exp(-\sqrt{\mathbf{s}^\top\mathbf{C}\mathbf{s}}). \tag{25}$$

Under the above conditions, the process $\mathbf{x}_t$ is multivariate Cauchy distributed with scale matrix $\mathbf{C}$. Consequently, for $0 < \alpha < \infty$, Example 2 gives the Rényi entropy of process $\mathbf{x}_t$ in the following form

$$R_\alpha(\mathbf{X}) = \begin{cases} \log\left(\frac{\alpha^{-\frac{d}{1-\alpha}}\sqrt{\det(4\pi\mathbf{C})}\,\Gamma(\frac{d+1}{2})}{\sqrt{\pi}}\right), & if \quad \alpha \neq 1; \\ \log\left(\frac{\sqrt{\det(4e^2\pi\mathbf{C})}\,\Gamma(\frac{d+1}{2})}{\sqrt{\pi}}\right), & if \quad \alpha = 1. \end{cases} \tag{26}$$

In practice, summations given in (19) and (22) can be evaluated from $j = 0$ to $m$, with $m \gg 1$ and $m < \infty$. This is because $\det(\mathbf{\Sigma}) \to 0$ when $m \to \infty$. Computational procedure for Rényi entropy of MCARMA systems based on multivariate Gaussian and Cauchy white noises is summarized in the steps described in Algorithm 5.2. From this procedure, it can be observed that computational complexity for Rényi entropy depends on estimation of $\mathbf{M}_l$ and $\mathbf{M}_l^*$, $l = 1, \dots, m$, for covariance matrix of MCARMA process. Therefore, Algorithm 5.2 has a computational complexity of $O(m^d)$, i.e., increments of steps of iteration $m$, increases the computational cost. Moreover, dimension $d$ of MCARMA system is often small (see Sections 7 and 8). Some disadvantages of Algorithm 5.2 are discussed in Section 9.

Rényi entropy computation of MCARMA system.

$d \times d$-matrices $\mathbf{S}_\mathbf{u}$, $\mathbf{S}_\mathbf{w}$, $\mathbf{A}_i$, $\mathbf{B}_j$, $\mathbf{D}_k$, $i = 1, \dots, n$, $j = 0, \dots, r$, $k = 0, \dots, q$. $m \gg 1$ and $\alpha > 0$, $\alpha \neq 1$. Set $\mathbf{M}_0 = \mathbf{M}_0^* = \mathbf{I}_d$. Compute $\mathbf{M}_l$ and $\mathbf{M}_l^*$ using (14) and (15), for $l = 1, \dots, m$. If $\mathbf{U} \sim N_d(\mathbf{0}, \mathbf{S}_\mathbf{u})$ and $\mathbf{W} \sim N_d(\mathbf{0}, \mathbf{S}_\mathbf{w})$, compute $\mathbf{\Sigma}$ and $R_\alpha(\mathbf{X})$ using (19) and (20). If $\mathbf{U} \sim C_d(\mathbf{0}, \mathbf{S}_\mathbf{u})$ and $\mathbf{W} \sim C_d(\mathbf{0}, \mathbf{S}_\mathbf{w})$, compute $\mathbf{C}$ and $R_\alpha(\mathbf{X})$ using (24) and (26).

### 5.3 Multivariate Laplace distribution

In this case, we take both the control and white noise stationary processes distributed multivariate Laplace with zero mean, $d \times d$ positive definite covariance matrices $\mathbf{S}_\mathbf{u}$ and $\mathbf{S}_\mathbf{w}$, respectively; thus both processes are respectively denoted as $L_d(\mathbf{0}, \mathbf{S}_\mathbf{u})$ and $L_d(\mathbf{0}, \mathbf{S}_\mathbf{w})$. The control and white noise processes satisfy the independence assumptions. The characteristic function of any random vector $\mathbf{y} \in \mathbb{R}^d$ which distributed multivariate Laplace with covariance matrix $\mathbf{\Sigma}_\mathbf{y}$ is in the following form [Fragiadakis and Meintanis(2011)]:

$$\boldsymbol{\varphi}_\mathbf{y}(\mathbf{s}) = \frac{1}{1+\frac{1}{2}\mathbf{s}^\top \mathbf{\Sigma}_\mathbf{y} \mathbf{s}}, \tag{27}$$

Therefore, the characteristic function of $\{\mathbf{x}_t\}$ can be written according to (17) as

$$\boldsymbol{\varphi}_{\mathbf{x}_t}(\mathbf{s}) = \prod_{j=0}^{\infty} \frac{1}{(1+\frac{1}{2}\mathbf{s}^\top \mathbf{K}_j \mathbf{s})(1+\frac{1}{2}\mathbf{s}^\top \mathbf{K}_j^* \mathbf{s})}, \tag{28}$$

where $\mathbf{K}_j = \mathbf{M}_j^\top \mathbf{S}_\mathbf{u} \mathbf{M}_j$ and $\mathbf{K}_j^* = \mathbf{M}^{*\top} \mathbf{S}_\mathbf{w} \mathbf{M}_j^*$, $j = 0,1,2,\dots.$ Expression (28) does not assign any traditional probability distribution of process $\mathbf{x}_t$, thus it is possible to use the uniqueness relation between the probability density and its characteristic function to get the distribution of $\mathbf{x}_t$. In fact, there is no analytical expression for the Rényi entropy of $\mathbf{x}_t$ related to control system (6) with component processes $\mathbf{u}_t$ and $\mathbf{w}_t$, therefore we consider the upper bound of Rényi entropy as follows.

## 6 Upper bound of Rényi entropy

According to determine the upper bound of Rényi entropy for the output process of stochastic control system in (6), we need to present expression of the covariance matrix of this process. In this section, we derive the formula of the covariance matrix to give us the upper bound of Rényi entropy for a MCARMA process which its distribution is difficult to determine.

From the representation of process $\mathbf{x}_t$ given in (11) and independence conditions, it can be seen that $\mathbb{E}[\mathbf{u}_t \mathbf{x}_s^\top] = 0$ and $\mathbb{E}[\mathbf{w}_t \mathbf{x}_s^\top] = 0$, for $t < s$. Therefore, for any $h >$ $\max\{r,q\}$, the autocovariance function $\mathbf{\Gamma}(h)=\mathbb{E}[\mathbf{x}_t\mathbf{x}_{t-h}^\top]$ of $\mathbf{x}_t$ [Brockwell and Davis(1991)] can be represented as follows

$$\mathbf{\Gamma}(h)=\sum_{i=1}^{p}\mathbf{A}_i\mathbf{\Gamma}(h-i). \tag{29}$$

Hence, latter relationship is used to compute the matrices $\mathbf{\Gamma}(h-i)$ for $h=p,p+1,\ldots,$ when $p<\max\{r,q\}$. We shall discuss next the case when $h<\max\{r,q\}$ by rewriting (6) as

$$\mathbf{X}_t=\mathbf{\Theta}\mathbf{X}_{t-1}+\mathbf{U}_t+\mathbf{W}_t, \tag{30}$$

where $\mathbf{X}_t=[\mathbf{x}_t^\top,\mathbf{x}_{t-1}^\top,\ldots,\mathbf{x}_{t-p+1}^\top,\mathbf{u}_t^\top,\mathbf{u}_{t-1}^\top,\ldots,\mathbf{u}_{t-r+1}^\top,\mathbf{w}_t^\top,\mathbf{w}_{t-1}^\top,\ldots,\mathbf{w}_{t-q+1}^\top]^\top$ is the $d(p+q+r)\times 1$ information vector and the parameter matrix is

$$\mathbf{\Theta}=\begin{pmatrix}\mathbf{\Theta}_{11} & \mathbf{\Theta}_{12}\\ \mathbf{\Theta}_{21} & \mathbf{\Theta}_{22}\end{pmatrix}, \tag{31}$$

where

$$\mathbf{\Theta}_{11}=\begin{pmatrix}\mathbf{A}_1 & \mathbf{A}_2 & \ldots & \mathbf{A}_p\\ \mathbf{I}_d & \mathbf{0} & \ldots & \mathbf{0}\\ \mathbf{0} & \ddots & \ldots & \vdots\\ \mathbf{0} & \mathbf{0} & \mathbf{I}_d & \mathbf{0}\end{pmatrix}_{dp\times dp},\ \mathbf{\Theta}_{12}=$$

$$\begin{pmatrix}\mathbf{B}_1 & \ldots & \mathbf{B}_r & \mathbf{D}_1 & \ldots & \mathbf{D}_q\\ \mathbf{0} & \ldots & \ldots & \ldots & \ldots & \mathbf{0}\\ \vdots & \vdots & \vdots & \vdots & \vdots & \vdots\\ \mathbf{0} & \ldots & \ldots & \ldots & \ldots & \mathbf{0}\end{pmatrix}_{dp\times d(r+q)},$$

$$\mathbf{\Theta}_{21}=\begin{pmatrix}\mathbf{0}\end{pmatrix}_{d(r+q)\times dp},\ \mathbf{\Theta}_{22}=\begin{pmatrix}\mathbf{R}_{11} & \mathbf{0}\\ \mathbf{0} & \mathbf{R}_{22}\end{pmatrix}_{d(r+q)\times d(r+q)},$$

$$\mathbf{R}_{11}=\begin{pmatrix}\mathbf{0} & \ldots & \ldots & \mathbf{0}\\ \mathbf{I}_d & \mathbf{0} & \ldots & \mathbf{0}\\ \mathbf{0} & \ddots & \ldots & \vdots\\ \mathbf{0} & \mathbf{0} & \mathbf{I}_d & \mathbf{0}\end{pmatrix}_{dr\times dr},\ \mathbf{R}_{22}=\begin{pmatrix}\mathbf{0} & \ldots & \ldots & \mathbf{0}\\ \mathbf{I}_d & \mathbf{0} & \ldots & \mathbf{0}\\ \mathbf{0} & \ddots & \ldots & \vdots\\ \mathbf{0} & \mathbf{0} & \mathbf{I}_d & \mathbf{0}\end{pmatrix}_{dq\times dq},$$

$$\mathbf{W}_t=\begin{pmatrix}\mathbf{w}_t & \mathbf{0} & \ldots & \mathbf{0} & \ldots & \mathbf{0} & \mathbf{w}_t & \mathbf{0} & \ldots & \mathbf{0}\end{pmatrix}_{1\times d(p+q+r)}^{\top},$$

$$\mathbf{U}_t=\begin{pmatrix}\mathbf{u}_t & \mathbf{0} & \ldots & \mathbf{u}_t & \mathbf{0} & \ldots & \mathbf{0} & \ldots & \mathbf{0}\end{pmatrix}_{1\times d(p+q+r)}^{\top}.$$

**Lemma 2** *The process $\boldsymbol{X}_t$ in the control system (30) is stable if and only if the process $\boldsymbol{x}_t$ in control system (6) is stable.*

**Proof.** For any complex number $z$, we have

$$\det(\mathbf{I}_{d(p+r+q)}-\mathbf{\Theta}z)=\det(\mathbf{I}_{dp}-\mathbf{\Theta}_{11}z)\det(\mathbf{I}_{d(r+q)}-\mathbf{\Theta}_{11}z).$$
$$=\det\left(\sum_{i=1}^{p}\mathbf{A}_i z^i\right).$$

We apply the representation (17) on (30) to get

$$\mathbf{X}_t = \sum_{j=0}^{\infty} \mathbf{\Theta}^j \mathbf{U}_{t-j} + \sum_{i=0}^{\infty} \mathbf{\Theta}^i \mathbf{W}_{t-i}. \quad (32)$$

Consider the following matrices

$$\tilde{\mathbf{I}} = \begin{pmatrix} \mathbf{I}_d & \mathbf{0} & \dots & \mathbf{0} \end{pmatrix}_{d\times d(p+r+q)},$$

$$\tilde{\mathbf{J}}_1 = \begin{pmatrix} \mathbf{I}_d & \mathbf{0} & \dots & \mathbf{I}_d & \mathbf{0} & \dots & \mathbf{0} & \dots & \mathbf{0} \end{pmatrix}^{\top}_{1\times d(p+q+r)},$$

$$\tilde{\mathbf{J}}_2 = \begin{pmatrix} \mathbf{I}_d & \mathbf{0} & \dots & \mathbf{0} & \dots & \mathbf{0} & \mathbf{I}_d & \mathbf{0} & \dots & \mathbf{0} \end{pmatrix}^{\top}_{1\times d(p+q+r)}.$$

Multiplying (32) by $\tilde{\mathbf{I}}$ and by some basic algebraic operations, we get

$$\begin{aligned}
\mathbf{X}_t &= \sum_{j=0}^{\infty} \tilde{\mathbf{I}}\mathbf{\Theta}^j \mathbf{U}_{t-j} + \sum_{i=0}^{\infty} \tilde{\mathbf{I}}\mathbf{\Theta}^i \mathbf{W}_{t-i} \\
&= \sum_{j=0}^{\infty} \tilde{\mathbf{I}}\mathbf{\Theta}^j \tilde{\mathbf{J}}_1 \tilde{\mathbf{I}} \mathbf{U}_{t-j} + \sum_{i=0}^{\infty} \tilde{\mathbf{I}}\mathbf{\Theta}^i \tilde{\mathbf{J}}_2 \tilde{\mathbf{I}} \mathbf{W}_{t-i} \\
&= \sum_{j=0}^{\infty} \tilde{\mathbf{I}}\mathbf{\Theta}^j \tilde{\mathbf{J}}_1 \mathbf{u}_{t-j} + \sum_{i=0}^{\infty} \tilde{\mathbf{I}}\mathbf{\Theta}^i \tilde{\mathbf{J}}_2 \mathbf{w}_{t-i}.
\end{aligned}$$

If we compare the above equation with equation in (11), then we get $\mathbf{M}_j = \tilde{\mathbf{I}}\mathbf{\Theta}^j\tilde{\mathbf{J}}_1$ and $\mathbf{M}_i^* = \tilde{\mathbf{I}}\mathbf{\Theta}^i\tilde{\mathbf{J}}_2$, for $j, i = 0,1, \ldots$. Latter suggest that control systems (6) and (30) are equivalent. This complete the proof of this Lemma. ☐

Multiplying by $\mathbf{X}_{t-h}^{\top}$, $h = 0,1,2, \ldots$, and taking the expectation for both sides of (30), we obtain

$$\mathbb{E}[\mathbf{X}_t\mathbf{X}_{t-h}] = \mathbf{\Theta}\mathbb{E}[\mathbf{X}_{t-1}\mathbf{X}_{t-h}] + \mathbb{E}[\mathbf{U}_t\mathbf{X}_{t-h}] + \mathbb{E}[\mathbf{W}_t\mathbf{X}_{t-h}].$$

If we consider the values $h = 0$ and 1 in latter equation, and use independent assumptions of system (6), then the covariance matrix $\tilde{\mathbf{\Gamma}}(0)$ of $\mathbf{X}_t$ is

$$\tilde{\mathbf{\Gamma}}(0) = \mathbf{\Theta}\tilde{\mathbf{\Gamma}}(1)^{\top} + \mathbf{S}_{\mathbf{U}} + \mathbf{S}_{\mathbf{W}}, \quad (33)$$

and, for $h = 1$,

$$\tilde{\mathbf{\Gamma}}(1) = \mathbf{\Theta}\tilde{\mathbf{\Gamma}}(0). \quad (34)$$

Substituting (34) in (33), we get

$$\tilde{\mathbf{\Gamma}}(0) = \mathbf{\Theta}\tilde{\mathbf{\Gamma}}(0)\mathbf{\Theta}^{\top} + \mathbf{S}_{\mathbf{U}} + \mathbf{S}_{\mathbf{W}}. \quad (35)$$

The vector operator for both sides of (35) gives

$$vec[\tilde{\mathbf{\Gamma}}(0)] = vec[\mathbf{\Theta}\tilde{\mathbf{\Gamma}}(0)\mathbf{\Theta}^{\top}] + vec[\mathbf{S}_{\mathbf{U}}] + vec[\mathbf{S}_{\mathbf{W}}].$$

Using the property of vector operator

$$vec[\mathbf{\Theta}\tilde{\mathbf{\Gamma}}(0)\mathbf{\Theta}^{\top}] = \mathbf{\Theta} \otimes \mathbf{\Theta}^{\top} vec[\tilde{\mathbf{\Gamma}}(0)],$$

where $\otimes$ is the Kronecker product, we get

$$(\mathbf{I}_{(d(p+r+q))^2} - \mathbf{\Theta} \otimes \mathbf{\Theta}^{\top}) vec[\tilde{\mathbf{\Gamma}}(0)] = vec[\mathbf{S}_{\mathbf{U}}] + vec[\mathbf{S}_{\mathbf{W}}]. \quad (36)$$

Since $\mathbf{x}_t$ is stable by Lemma 2, then the matrix $\mathbf{I}_{(d(p+r+q))^2} - \mathbf{\Theta} \otimes \mathbf{\Theta}^{\top}$ is invertible. Therefore

$$vec[\tilde{\mathbf{\Gamma}}(0)] = (\mathbf{I}_{(d(p+r+q))^2} - \mathbf{\Theta} \otimes \mathbf{\Theta}^{\top})^{-1} vec[\mathbf{S}_{\mathbf{U}}] + (\mathbf{I}_{(d(p+r+q))^2} - \mathbf{\Theta} \otimes \mathbf{\Theta}^{\top})^{-1} vec[\mathbf{S}_{\mathbf{W}}]. \quad (37)$$

After computing $vec[\tilde{\mathbf{\Gamma}}(0)]$ from (37), we can collect $\mathbf{\Gamma}(0), \mathbf{\Gamma}(1), \ldots, \mathbf{\Gamma}(p-1)$, from the representation

$$\tilde{\mathbf{\Gamma}}(0) = \begin{pmatrix} \mathbf{\Gamma}_{11} & \mathbf{\Gamma}_{12} \\ \mathbf{\Gamma}_{12}^{\top} & \mathbf{\Gamma}_{22} \end{pmatrix}, \quad (38)$$

where

$$\mathbf{\Gamma}_{11}=\begin{pmatrix}\mathbf{\Gamma}(0) & \mathbf{\Gamma}(1) & \dots & \mathbf{\Gamma}(p-1)\\ \mathbf{\Gamma}(-1) & \mathbf{\Gamma}(0) & \dots & \mathbf{\Gamma}(p-1)\\ \vdots & \vdots & \ddots & \vdots\\ \mathbf{\Gamma}(-p+1) & \mathbf{\Gamma}(-p+2) & \dots & \mathbf{\Gamma}(0)\end{pmatrix},$$

$$\mathbf{\Gamma}_{12}=\begin{pmatrix}k_{u0} & k_{u1} & \dots & k_{ur} & k_{w0} & k_{w1} & \dots & k_{wq}\end{pmatrix},$$

$$\mathbf{\Gamma}_{22}=\begin{pmatrix}\mathbf{S}_{11} & \mathbf{0}\\ \mathbf{0} & \mathbf{S}_{22}\end{pmatrix},$$

$$\mathbf{S}_{11}=\begin{pmatrix}\mathbf{S}_u & \mathbf{0} & \dots & \mathbf{0}\\ \mathbf{0} & \ddots & \dots & \vdots\\ \vdots & \vdots & \ddots & \vdots\\ \mathbf{0} & \mathbf{0} & \dots & \mathbf{S}_u\end{pmatrix},\mathbf{S}_{22}=\begin{pmatrix}\mathbf{S}_w & \mathbf{0} & \dots & \mathbf{0}\\ \mathbf{0} & \ddots & \dots & \vdots\\ \vdots & \vdots & \ddots & \vdots\\ \mathbf{0} & \mathbf{0} & \dots & \mathbf{S}_w\end{pmatrix},$$

$$k_{u0}=\begin{pmatrix}\mathbb{E}[\mathbf{x}_t\mathbf{u}_t^\top] & \mathbf{0} & \dots & \mathbf{0}\end{pmatrix}^\top,$$

$$k_{u1}=\begin{pmatrix}\mathbb{E}[\mathbf{x}_t\mathbf{u}_{t-1}^\top] & \mathbb{E}[\mathbf{x}_{t-1}\mathbf{u}_{t-1}^\top] & \mathbf{0} & \dots & \mathbf{0}\end{pmatrix}^\top,$$

$$\vdots$$

$$k_{ur}=\begin{pmatrix}\mathbb{E}[\mathbf{x}_t\mathbf{u}_{t-r+1}^\top] & \mathbb{E}[\mathbf{x}_{t-1}\mathbf{u}_{t-r+1}^\top] & \mathbb{E}[\mathbf{x}_{t-2}\mathbf{u}_{t-r+1}^\top] & \dots & \mathbb{E}[\mathbf{x}_{t-p}\mathbf{u}_{t-r+1}^\top]\end{pmatrix}^\top,$$

$$k_{w0}=\begin{pmatrix}\mathbb{E}[\mathbf{x}_t\mathbf{w}_t^\top] & \mathbf{0} & \dots & \mathbf{0}\end{pmatrix}^\top,$$

$$k_{w1}=\begin{pmatrix}\mathbb{E}[\mathbf{x}_t\mathbf{w}_{t-1}^\top] & \mathbb{E}[\mathbf{x}_{t-1}\mathbf{w}_{t-1}^\top] & \mathbf{0} & \dots & \mathbf{0}\end{pmatrix}^\top,$$

$$\vdots$$

$$k_{wr}=\begin{pmatrix}\mathbb{E}[\mathbf{x}_t\mathbf{w}_{t-r+1}^\top] & \mathbb{E}[\mathbf{x}_{t-1}\mathbf{w}_{t-r+1}^\top] & \mathbb{E}[\mathbf{x}_{t-2}\mathbf{w}_{t-r+1}^\top] & \dots & \mathbb{E}[\mathbf{x}_{t-p}\mathbf{w}_{t-r+1}^\top]\end{pmatrix}^\top.$$

Given $\mathbf{\Gamma}(0)$, the following Corollary is straightforward from Lemma 1.

**Corollary 1** *The upper bound for the Rényi entropy of MCARMA process $\boldsymbol{x}_t$ defined in (6) is*

$$R_\alpha^{Upper}(\mathbf{X})=\mathcal{C}_d(\alpha)+\frac{1}{2}\log[\det(\mathbf{\Gamma}(0))], \tag{39}$$

where $\mathcal{C}_d(\alpha)$ is defined in (3).

# 7 Numerical examples

In this section, the examples are given to show the effectiveness of characteristic function (17) and upper bound (39). In particular, consider the MCARMA system of dimension $d = 3$ as follows

$$\mathbf{x}_t = \mathbf{A}_1\mathbf{x}_{t-1} + \mathbf{u}_t + \mathbf{B}_1\mathbf{u}_{t-1} + \mathbf{w}_t + \mathbf{D}_1\mathbf{w}_{t-1}, \tag{40}$$

with $\mathbf{x}_0 = \mathbf{0}$, whose both processes $\mathbf{u}_t$ and $\mathbf{w}_t$ satisfy the independence assumptions and have distributions depending on the following cases.

### 7.1 Multivariate Gaussian case

Consider process (40), where $\mathbf{u}_t \in \mathbb{R}^3$ is a zero mean stationary Gaussian process with covariance matrix $\mathbf{S}_\mathbf{u}$ which represents the system input vector, $\mathbf{w}_t \in \mathbb{R}^3$ is the stationary Gaussian white noise process with mean zero and covariance matrix $\mathbf{S}_\mathbf{w}$. Consider the following parameters

$$\mathbf{A}_1 = \begin{pmatrix} 0.5 & 0 & 0 \\ 0.1 & 0.1 & 0.3 \\ 0 & 0.2 & 0.3 \end{pmatrix}, \mathbf{B}_1 = \begin{pmatrix} 0.3 & 0 & 0 \\ 0 & 0.1 & 0.2 \\ 0 & 0.2 & 0.3 \end{pmatrix},$$

$$\mathbf{D}_1 = \mathbf{I}_3,$$

$$\mathbf{S}_\mathbf{u} = \begin{pmatrix} 2.25 & 0 & 0 \\ 0 & 1 & 0.5 \\ 0 & 0.5 & 0.74 \end{pmatrix}, \mathbf{S}_\mathbf{w} = \begin{pmatrix} 0.25 & 0 & 0 \\ 0 & 1 & 0 \\ 0 & 0 & 0.5 \end{pmatrix}.$$

The characteristic polynomial determinant of this process is $\det(\mathbf{I}_3 - \mathbf{A}_1 z) = (1 - 0.5z)(1 - 0.4z - 0.03z^2)$. Therefore, the roots of these polynomials are 2, 2.153, $-15.486$. Clearly, the length of all roots greater than 1 and the output process $\mathbf{x}_t$ is stable. The matrices $\mathbf{M}_i$ and $\mathbf{M}_i^*$, $i = 0,1,\dots$, are computed from (14) and (15) as

$$\mathbf{M}_0 = \begin{pmatrix} 0.4 & 0 & 0 \\ 0.09 & 0.14 & 0.23 \\ 0.02 & 0.16 & 0.28 \end{pmatrix}, \mathbf{M}_0^* = \begin{pmatrix} 0.75 & 0 & 0 \\ 0.16 & 0.17 & 0.42 \\ 0.02 & 0.28 & 0.45 \end{pmatrix},$$

$$\mathbf{M}_1 = \begin{pmatrix} 0.2 & 0 & 0 \\ 0.055 & 0.062 & 0.107 \\ 0.024 & 0.076 & 0.13 \end{pmatrix}, \mathbf{M}_1^* = \begin{pmatrix} 0.375 & 0 & 0 \\ 0.097 & 0.101 & 0.177 \\ 0.038 & 0.118 & 0.219 \end{pmatrix},$$

$$\mathbf{M}_2 = \begin{pmatrix} 0.1 & 0 & 0 \\ 0.0327 & 0.029 & 0.0497 \\ 0.0182 & 0.0352 & 0.0604 \end{pmatrix}, \mathbf{M}_2^* = \begin{pmatrix} 0.1875 & 0 & 0 \\ 0.0586 & 0.0455 & 0.0834 \\ 0.0308 & 0.0556 & 0.1011 \end{pmatrix},$$

$$\mathbf{M}_3 = \begin{pmatrix} 0.05 & 0 & 0 \\ 0.0187 & 0.0135 & 0.0231 \\ 0.0120 & 0.0164 & 0.0281 \end{pmatrix}, \mathbf{M}_3^* = \begin{pmatrix} 0.0938 & 0 & 0 \\ 0.0339 & 0.0212 & 0.0387 \\ 0.0210 & 0.0258 & 0.0470 \end{pmatrix},$$

$$\vdots$$

Latter matrices are used to obtain $\mathbf{\Sigma}$ of (19) and $\det(\mathbf{\Sigma}) \approx 3.943$. Then, Equation (20) gives the values for Rényi entropy of process $\mathbf{x}_t$ as follow

$$R_\alpha(\mathbf{X}) \approx \begin{cases} 3.443 - \frac{3}{2}\frac{\log(\alpha)}{1-\alpha}, & if \quad \alpha \neq 1; \\ 4.943, & if \quad \alpha \to 1. \end{cases} \tag{41}$$

Rényi entropy is upper bounded by equation (39) as follows. To compute $vec(\tilde{\mathbf{\Gamma}}(0))$ in (37), the matrices $\mathbf{\Gamma}(0), \mathbf{\Gamma}(1), \ldots, \mathbf{\Gamma}(p-1)$ are collected from representation (38). Then

$$\mathbf{\Gamma}(0) = \begin{pmatrix} 5.1700 & 0.3765 & 0.0443 \\ 0.3765 & 0.9241 & 1.3560 \\ 0.0443 & 1.3560 & 2.8917 \end{pmatrix}.$$

Therefore, from (39), the upper bound is

$$R_\alpha^{Upper}(\mathbf{X}) \approx C_3(\alpha) + 3.943. \tag{42}$$

## 7.2 Multivariate Cauchy case

Consider ARMA control system (40) with the following assumptions: $\mathbf{u}_t \in \mathbb{R}^3$ is a zero mean stationary Cauchy process with scale matrix $\mathbf{S}_\mathbf{u}$ which it represents the system input vector, $\mathbf{w}_t \in \mathbb{R}^3$ is the stationary Cauchy white noise process with mean zero and scale matrix $\mathbf{S}_\mathbf{w}$. Consider the following parameters:

$$\mathbf{A}_1 = 0.5\mathbf{I}_3, \mathbf{B}_1 = 0.3\mathbf{I}_3,$$
$$\mathbf{D}_1 = \mathbf{I}_3,$$
$$\mathbf{S}_\mathbf{u} = \begin{pmatrix} 0.25 & 0 & 0 \\ 0 & 1 & 0 \\ 0 & 0 & 0.5 \end{pmatrix}, \mathbf{S}_\mathbf{w} = \begin{pmatrix} 0.25 & 0 & 0 \\ 0 & 1 & 0 \\ 0 & 0 & 0.5 \end{pmatrix}.$$

The characteristic polynomial determinant of this process is $\mathrm{d}et(\mathbf{I}_3 - \mathbf{A}_1 z) = (1 - 0.5z)^3$. Therefore, all roots of this polynomials are 2. Then, the length of all roots greater than 1 and the output process $\mathbf{x}_t$ is stable. The matrices $\mathbf{M}_i$ and $\mathbf{M}_i^*$, $i = 0,1, \ldots$, are computed from (14) and (15) as

$$\mathbf{M}_0 = 0.4\mathbf{I}_3, \mathbf{M}_0^* = 0.75\mathbf{I}_3,$$
$$\mathbf{M}_1 = 0.2\mathbf{I}_3, \mathbf{M}_1^* = 0.375\mathbf{I}_3,$$
$$\mathbf{M}_2 = 0.1\mathbf{I}_3, \mathbf{M}_2^* = 0.1875\mathbf{I}_3,$$
$$\mathbf{M}_3 = 0.05\mathbf{I}_3, \mathbf{M}_3^* = 0.0938\mathbf{I}_3,$$
$$\vdots$$
$$\mathbf{L}_\mathbf{u} = \mathbf{L}_\mathbf{w} = \begin{pmatrix} 0.75 & 0 & 0 \\ 0.3162 & 0 & 0 \\ 0 & 0 & 0.7071 \end{pmatrix},$$
$$\mathbf{K}_0 = 0.4\mathbf{I}_3\mathbf{L}_\mathbf{u}, \mathbf{K}_0^* = 0.75\mathbf{I}_3\mathbf{L}_\mathbf{w},$$
$$\mathbf{K}_1 = 0.2\mathbf{I}_3\mathbf{L}_\mathbf{u}, \mathbf{K}_1^* = 0.375\mathbf{I}_3\mathbf{L}_\mathbf{w},$$
$$\mathbf{K}_2 = 0.1\mathbf{I}_3\mathbf{L}_\mathbf{u}, \mathbf{K}_2^* = 0.1875\mathbf{I}_3\mathbf{L}_\mathbf{w},$$
$$\mathbf{K}_3 = 0.05\mathbf{I}_3\mathbf{L}_\mathbf{u}, \mathbf{K}_3^* = 0.0938\mathbf{I}_3\mathbf{L}_\mathbf{w},$$
$$\mathbf{K}_4 = 0.025\mathbf{I}_3\mathbf{L}_\mathbf{u}, \mathbf{K}_4^* = 0.0469\mathbf{I}_3\mathbf{L}_\mathbf{w},$$
$$\vdots$$

It is clear that proportional condition of matrices $\mathbf{K}_j\mathbf{K}_j^\top$, $\mathbf{K}_j^*\mathbf{K}_j^{*\top}$, $j = 0,1,2, \ldots$, is hold. Then,

$\mathbf{C} = a^2\mathbf{S_u}$ and $a = (0.4 + 0.2 + 0.1 + 0.05 + 0.025 + \cdots + 0.75 + 0.375 + 0.1875 + 0.0938 + \cdots)$. Then, $\det(\mathbf{C}) \approx 1.850397$ and, from (26), the values of Rényi entropy of process $\mathbf{x}_t$ are

$$R_\alpha(\mathbf{X}) \approx \begin{cases} 3.5319 - \frac{3}{1-\alpha}\log(\alpha), & if \quad \alpha \neq 1; \\ 6.5319, & if \quad \alpha \to 1. \end{cases} \tag{43}$$

### 7.3 Multivariate Laplace case

In this case, we consider the processes $\mathbf{u}_t$ and $\mathbf{w}_t$ with multivariate Laplace distribution. From (39), the upper bound for Rényi of $\mathbf{x}_t$ is given by (42).

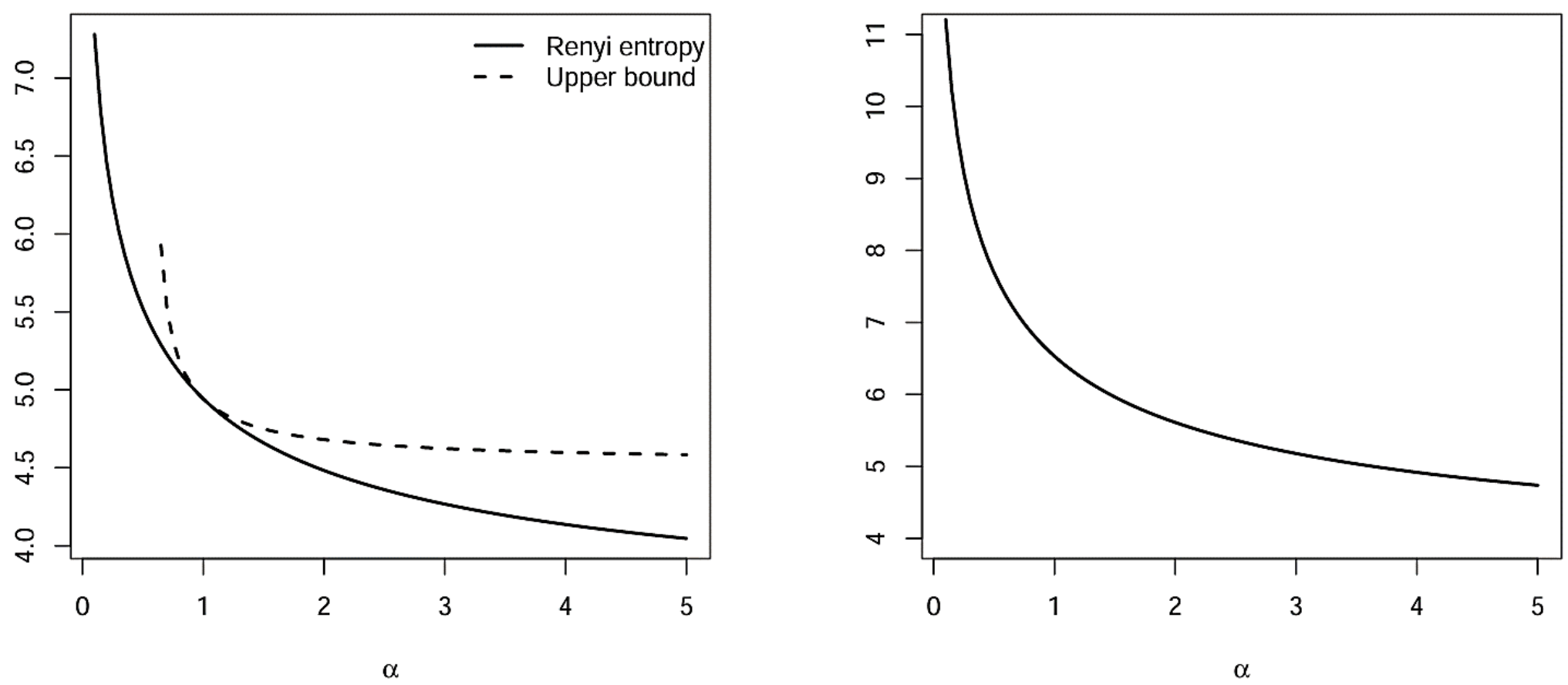

Figure 1: Left: Rényi entropy and their respective upper bound for multivariate Gaussian example. Right: Rényi entropy for multivariate Cauchy example.

Left plot of Figure 1 illustrates the performance of multivariate Gaussian Rényi entropy for $0 < \alpha \leq 5$ and their respective upper bound. Note that some values of the upper bound were not plotted because, from (3), $C_3(\alpha)$ is not available for $\alpha \leq d/(d+2) = 3/5$. Nevertheless, for other values of $\alpha$, upper bound is relatively close to Rényi entropy. These upper bound values are also considered for Laplace case. Moreover, Rényi entropy has a decreasing behavior when $\alpha$ increases. Right plot shows a similar behavior for multivariate Cauchy case; however, these values are larger than multivariate Gaussian ones.

## 8 Applications

### 8.1 Paper-making process

Consider the paper-making process considered in [Liu et al.(2018)] modeled by a multi-input and multi-output stochastic system MCARMA system. The paper-making process is summarized as follows: first, the paper fiber materials are pumped into a flow box to be diluted with water. After discharging from the box, the second step is to drain water from fiber materials. Third, the pressing machine formed the paper. Fourth, using steam heating the remaining water is removed from the paper. Assume that the system from can be modeled by the following multivariate model:

$$\mathbf{x}_t + \mathbf{A}_1\mathbf{x}_{t-1} + \mathbf{A}_2\mathbf{x}_{t-2} = \mathbf{B}_1\mathbf{u}_{t-1} + \mathbf{B}_2\mathbf{u}_{t-2} + \mathbf{D}_1\mathbf{w}_{t-1}, \qquad (44)$$

with $\mathbf{x}_0 = \mathbf{0}$, $\mathbf{D}_1 = \mathbf{I}_2$, and both processes $\mathbf{u}_{t-1}$ and $\mathbf{w}_{t-1}$ satisfy the independence assumptions and have multivariate Gaussian distributions. Here, $\mathbf{u}_t = [u_1(t), u_2(t)]^\top$ involves the amount of the materials $u_1(t)$ and the steam pressure $u_2(t)$, respectively; and $\mathbf{x}_t = [x_1(t), x_2(t)]^\top$ involves the moisture $x_1(t)$ and the basis weight $x_2(t)$. [Liu et al.(2018)] considered particular values for components of matrices $\mathbf{A}_1$, $\mathbf{A}_2$, $\mathbf{B}_1$ and $\mathbf{B}_2$ in system (44), thus the corresponding matrices are:

$$\mathbf{A}_1 = \begin{pmatrix} -0.62 & 0.45 \\ -0.58 & -0.32 \end{pmatrix}, \mathbf{B}_1 = \begin{pmatrix} -0.86 & 0.16 \\ 0.68 & -0.44 \end{pmatrix},$$

$$\mathbf{A}_2 = \begin{pmatrix} -0.33 & 0.29 \\ -0.24 & 0.48 \end{pmatrix}, \mathbf{B}_2 = \begin{pmatrix} -0.84 & -0.65 \\ 0.55 & 0.61 \end{pmatrix}.$$

Figure 2 illustrates the performance of multivariate Gaussian Rényi entropy for $0 < \alpha \leq 5$ for process (44). As multivariate Gaussian case, the Rényi entropy has a decreasing behavior when $\alpha$ increases. However, when $\alpha = 1$ (Shannon entropy case), the lowest information value is found.

Rényi entropy measure could be useful to compare information of different systems obtained from different settings of matrices $\mathbf{A}_1$, $\mathbf{A}_2$, $\mathbf{B}_1$ and $\mathbf{B}_2$ in (44). In this sense, [Liu et al.(2018)] proposed the standard multivariable generalized extended stochastic gradient (M-GESG) algorithm. Therefore proposed Rényi entropy is suitable to measure M-GESG's performance; however, given that authors obtained stable and small parameter estimation errors when $t$ increases, this step is skipped in this illustration.

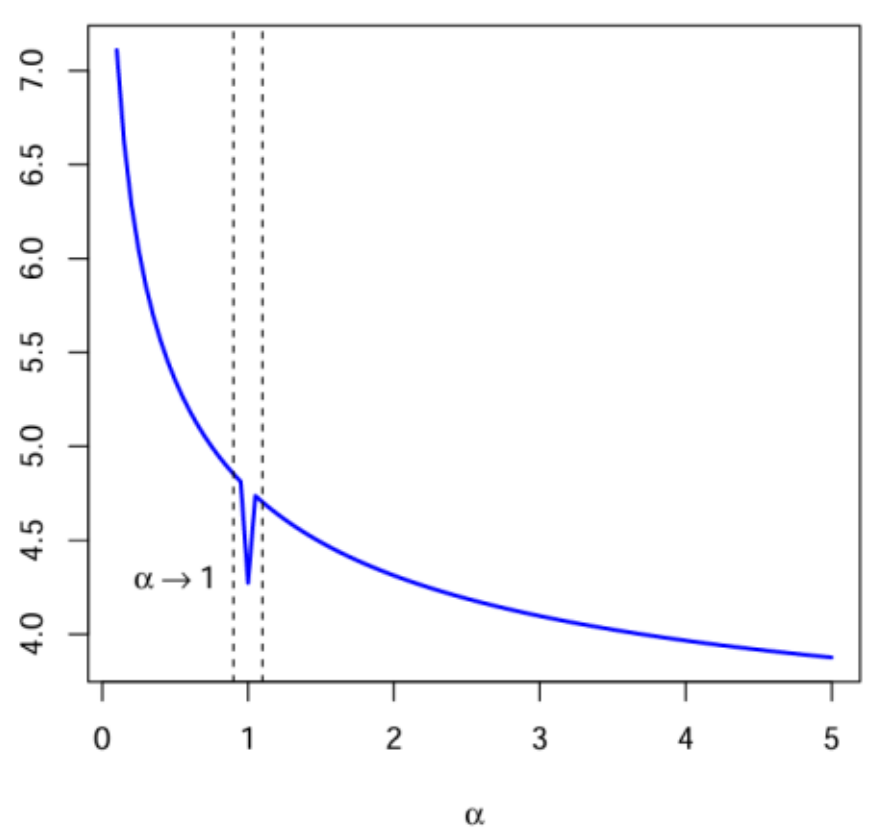

Figure 2: Rényi entropy based on multivariate Gaussian white noise for paper-making process.

## 8.2 Electric circuit model

Consider the MCARMA system for an electric circuit model with two-input voltage and a single-output voltage measured from a capacitor proposed by [Hag ElAmin(2020)]. From the circuit system, the author considered the following model:

$$\mathbf{x}_t + \mathbf{A}\mathbf{x}_{t-1} = \mathbf{B}\mathbf{u}_{t-1} + \mathbf{D}\mathbf{w}_{t-1}, \tag{45}$$

with $\mathbf{x}_0 = \mathbf{0}$. Here, an additive noise $\mathbf{w}_{t-1}$ is included with $\mathbf{D} = \mathbf{I}_2$, and both processes $\mathbf{u}_{t-1}$ and $\mathbf{w}_{t-1}$ satisfy the independence assumptions and have multivariate Gaussian distributions. The process $\mathbf{x}_t = [x_1(t), x_2(t)]^\top$ involves the inductor current $x_1(t)$ as the first state and the capacitor voltage $x_2(t)$ as the second state. The corresponding matrices are:

$$\mathbf{A} = \begin{pmatrix} \frac{R_2}{L} & \frac{1}{L} \\ -\frac{1}{C} & \frac{1}{C}\left[\frac{1}{R_1} + \frac{1}{R_3}\right] \end{pmatrix}, \mathbf{B} = \begin{pmatrix} 0 & \frac{1}{C} \\ \frac{1}{R_1 C} & 0 \end{pmatrix},$$

where $R_i$, $i = 1,2,3$, are the resistors, $L$ is the inductor, and $C$ is the capacitor. In particular, [Hag ElAmin(2020)] choose $R_1 = 2$ $(\Omega)$, $R_2 = 2.95$ $(\Omega)$, $R_3 = 1.15$ $(\Omega)$, $L = 2.6788$ $(H)$, and $C = 2.2$ $(F)$.

Figure 3 illustrates the uncertainty quantification of process (45) using the multivariate Gaussian Rényi entropy. In this case, some variations of parameter $L$ were made. Rényi entropy has a decreasing behavior when $\alpha$ increases and the lowest information value is found when $\alpha = 1$ (Shannon entropy case). However, entropy decreases when inductor $L$ increase because, approximation of matrix $\mathbf{\Sigma}$ depends on components of matrices $\mathbf{A}$ and $\mathbf{B}$. Therefore, Rényi entropy measures is affected mainly for large values of inductor $L$.

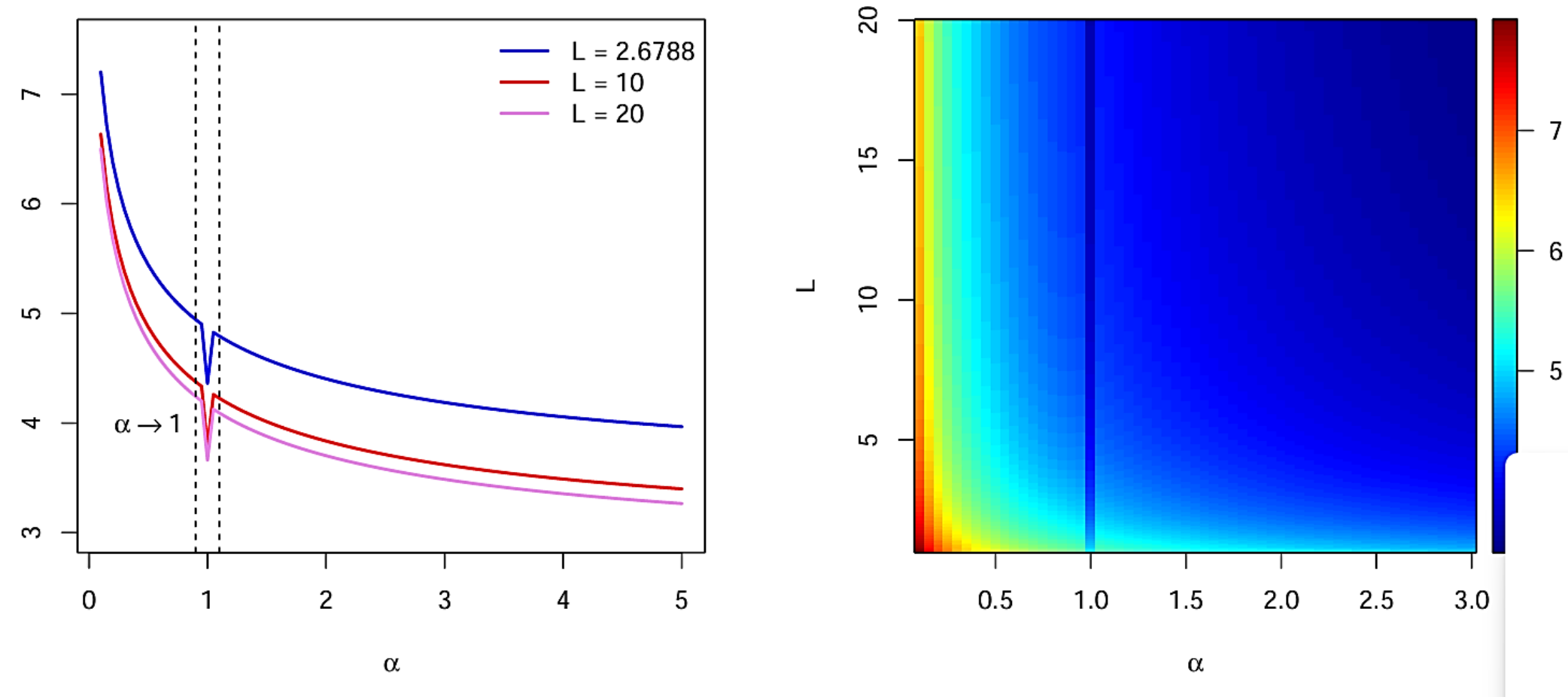

Figure 3: Multivariate Gaussian Rényi entropy for electric circuit model with some specific (left) and varying (right) $L$ parameters.

# 9 Conclusions

In this paper, we have derived the formula of the characteristic function of MCARMA control process, and then we used this formula as a tool to find the expression of Rényi entropy. As special case, three models (Gaussian, Cauchy and Laplace) are discussed. Also the covariance matrix and upper bound of entropy are computed. Lemma 2 analyzes the necessary conditions for stability of MCARMA process, which is required for the upper bound computation for the Rényi entropy. In addition, computational cost for Rényi entropy computation of MCARMA system presented in this paper depends mainly on number of iterations $m$ (until convergence of $det(\mathbf{\Sigma})$ and $det(\mathbf{C})$ for multivariate Gaussian and Cauchy cases, respectively) and MCARMA system dimension $d$. Thus large values of $m$ and $d$ increases complexity of Rényi entropy computation.

Classical methods for uncertainty/information quantification are focused on multivariate autoregressive-type models [Gibson(2018), Contreras-Reyes(2022a), Contreras-Reyes(2024a)], where the output data parameters of system identification are missing. Therefore, this study is pioneer in uncertainty/information quantification of MCARMA systems, where both input and output parameters are considered. In particular, in Section 8, the Rényi entropy quantified the information provided by a multi-input and multi-output paper-making stochastic system and an electric circuit model with two-input voltage and a single-output voltage system. Moreover, the Rényi entropy summarized the parameter set information of these industrial processes.

The proposed Rényi information measure of MCARMA system can be further considered for data analysis in various systems such as speech processing, biological signal processing, electroencephalogram classification, electrocardiogram classification, geophysical exploration, communication systems, and environmental issues. Moreover, Rényi entropy is available to quantify information of MCARMA system on each iteration of least square estimation methods [Bao et al.(2011), Liu et al.(2018)] and determine model performance.

## Declaration of Competing Interest

The authors declare that they have no known competing financial interests or conflict of interest that could have appeared to influence the work reported in this paper.

## Data availability statement

Data sharing not applicable – no new data generated.

## Acknowledgements

The author thanks the editor and three anonymous referees for their helpful comments and suggestions. R codes used in this paper are available upon request from the corresponding author.